\documentclass[12pt]{article}
\usepackage{graphicx}
\usepackage{amsmath}
\usepackage{amsthm}
\usepackage{amsfonts}
\usepackage{amssymb, amscd}

\newtheorem{thm}{Theorem}[section]
\newtheorem{cor}[thm]{Corollary}
\newtheorem{lem}[thm]{Lemma}
\newtheorem{prop}[thm]{Proposition}

\theoremstyle{definition}
\newtheorem{defn}[thm]{Definition}
\theoremstyle{remark}
\newtheorem{rem}[thm]{Remark}
\numberwithin{equation}{section}

\newcommand{\K}{\mathbb K}

\newcommand{\A}{\mathcal{A}}

\begin{document}

\begin{center}
{\large \bf Hom-algebra structures}%
\bigskip

{\bf Abdenacer Makhlouf} \\
{\rm \small Universit\'{e} de Haute Alsace,  \\
Laboratoire de Math\'{e}matiques, Informatique et Applications, \\
4, rue des Fr\`{e}res Lumi\`{e}re F-68093 Mulhouse, France} \\
Abdenacer.Makhlouf@uha.fr \\
\medskip

{ \bf Sergei Silvestrov} \\
{\rm \small Centre for Mathematical Sciences, Lund University, \\
Box 118, SE-221 00 Lund, Sweden} \\
sergei.silvestrov@math.lth.se  \end{center}


\bigskip

\begin{center}  {\bf Abstract} \end{center}

A Hom-algebra structure is a multiplication on a
vector space where the structure is twisted by a
homomorphism. The structure of Hom-Lie algebra
was introduced by Hartwig, Larsson and Silvestrov
in \cite{HLS} and extended by Larsson and
Silvestrov to quasi-hom Lie and quasi-Lie
algebras in \cite{LS1,LS2}. In this paper we
introduce and study Hom-associative, Hom-Leibniz,
and Hom-Lie admissible algebraic structures which
generalize the well known associative, Leibniz
and Lie admissible algebras. Also, we
characterize the flexible Hom-algebras in this
case. We also explain some connections  between
Hom-Lie algebras and Santilli's isotopies of
associative and Lie algebras.

\noindent {\bf AMS MSC 2000:} 17A30,16Y99,17A01,17A20,17D25 \\
{\bf Keywords:} Hom-Lie algebra, Hom-Associative
algebra, Hom-Leibniz algebra, Hom-Lie admissible
algebra, flexible algebra,
classification \\


\newpage

\section*{Introduction}

In \cite{HLS,LS1,LS2}, the class of quasi-Lie algebras and
subclasses of quasi-hom-Lie algebras and Hom-Lie algebras has been
introduced. These classes of algebras are tailored in a way suitable
for simultaneous treatment of the Lie algebras,  Lie superalgebras,
 color Lie algebras and  deformations arising in connection
with twisted, discretized or deformed derivatives and corresponding
generalizations, discretizations and deformations of vector fields
and differential calculus. It has been shown in
\cite{HLS,LS1,LS2,LS3} that the class of quasi-hom-Lie algebras
contains as a subclass on one hand the color Lie algebras and in
particular Lie superalgebras and Lie algebras, and on another hand
various known and new single and multi-parameter families of
algebras obtained using twisted derivations and constituting
deformations and quasi-deformations of universal enveloping algebras
of Lie and color Lie algebras and of algebras of vector-fields. The
main feature of quasi-Lie algebras, quasi-hom-Lie algebras and
Hom-Lie algebras is that the skew-symmetry and the Jacobi identity
are twisted by several deforming twisting maps and also in quasi-Lie
and quasi-hom-Lie algebras the Jacobi identity in general contains 6
twisted triple bracket terms.

In this paper, we provide a different way for
constructing Hom-Lie algebras by extending the
fundamental construction of Lie algebras from
associative algebras via commutator bracket
multiplication. To this end we define the notion
of Hom-associative algebras generalizing
associative algebras to a situation where
associativity law is twisted, and show that the
commutator product defined using the
multiplication in a Hom-associative algebra leads
naturally to Hom-Lie algebras. We introduce also
Hom-Lie admissible algebras and more general
$G$-Hom-associative algebras with subclasses of
Hom-Vinberg and pre-Hom-Lie algebras,
generalizing to the twisted situation Lie
admissible algebras, $G$-associative algebras,
Vinberg and pre-Lie algebras respectively, and
show that for these classes of algebras the
operation of taking commutator leads to Hom-Lie
algebras as well. We construct all the twistings
so that the brackets $[X_1,X_2]=2 X_2, \
[X_1,X_3]=-2 X_3, \ [X_2,X_3]=X_1$ determine a
three dimensional Hom-Lie algebra, generalizing
$sl(2)$. Finally, we provide for a subclass of
twistings, a list of all three-dimensional
Hom-Lie algebras. For some values of structure
constants, this list  contains all
three-dimensional Lie algebras. The families of
Hom-Lie algebras in these list can be viewed as
deformations of Lie algebras into a class of
Hom-Lie algebras.

The paper is organized as follows. In Section 1, we summarize the definitions of  Hom-associative,
 Hom-Leibniz algebra and Hom-Lie algebra. In Section 2, we extend the classical concept of
 Lie admissible algebra to Hom-Lie admissible algebra. We also review and show the connections to Santilli's
  isotopies of associative and Lie algebras. In Section 3, we explore some general classes of Hom-Lie admissible algebras.
  For any subgroup $G$ of the permutation group $\mathcal{S}_3$, we introduce
  the $G$-Hom-associative algebras, show that they are Hom-Lie admissible and describe all these classes.
  Section 4 is devoted to some properties of
  flexible Hom-algebras. We extend the notion of flexibility to Hom-algebras and characterize flexible Hom-algebras using
  Jordan and Lie parts of the multiplication. In Section 5, we consider the algebraic
  varieties of finite-dimensional Hom-associative and Hom-Lie algebras. We also provide a characterization of all Hom-Lie
algebras of $sl(2)$ type and construct all the $3$-dimensional
Hom-Lie algebras for a class of homomorphisms.

\noindent {\bf Acknowledgments.} This work was
partially supported by the Crafoord foundation,
The Royal Physiographic Society in Lund, The
Swedish Foundation of International Cooperation
in Research and Higher Education (STINT), The
Swedish Royal Academy of Sciences and the
Liegrits network.

\section{Hom-associative, Hom-Leibniz algebras and Hom-Lie algebras} Let
$\K$ be an algebraically closed field of
characteristic $0$ and $V$ be a linear space over
$\K$. Let $\alpha$ be a homomorphism of $V$.
\begin{defn}
A Hom-associative algebra is a triple $(V, \mu, \alpha)$ consisting
of a linear space $V$, a bilinear map $\mu: V\times V \rightarrow V$
and a  homomorphism $\alpha: V \rightarrow V$ satisfying
$$
\mu(\alpha(x),\mu (y,z))=\mu (\mu (x,y),\alpha (z))
$$
\end{defn}

A class of quasi-Leibniz algebras was introduced
in \cite{LS2} in connection to general quasi-Lie
algebras following the standard Loday's
conventions for Leibniz algebras (i.e. right
Loday algebras). In the context of the subclass
of Hom-Lie algebras one gets a class of
Hom-Leibniz algebras.

\begin{defn}
A Hom-Leibniz algebra is a triple $(V, [\cdot, \cdot], \alpha)$
consisting of a linear space $V$, bilinear map $[\cdot, \cdot]:
V\times V \rightarrow V$ and a homomorphism $\alpha: V \rightarrow
V$  satisfying
\begin{equation} \label{Leibnizalgident}
 [[x,y],\alpha(z)]=[[x,z],\alpha (y)]+[\alpha(x),[y,z]]
\end{equation}
\end{defn}

In terms of the (right) adjoint homomorphisms $Ad_Y: V\rightarrow V$
defined by $Ad_{Y}(X)=[X,Y]$, the identity \eqref{Leibnizalgident}
can be written as
\begin{equation} \label{LeibnizalgidentAdDeriv}
Ad_{\alpha(z)}([x,y]) = [Ad_{z}(x),\alpha(y)] +[\alpha(x),Ad_{z}(y)]
\end{equation}
or in pure operator form
\begin{equation} \label{LeibnizalgidentAdOper}
Ad_{\alpha(z)} \circ Ad_{y} = Ad_{\alpha(y)} \circ Ad_{z} +
Ad_{Ad_{z}(y)} \circ \alpha
\end{equation}


\subsection{Hom-Lie algebras}
The Hom-Lie algebras were initially introduced by Hartwig, Larson
and Silvestrov in \cite{HLS} motivated initially by examples of
deformed Lie algebras coming from twisted discretizations of vector
fields.

\begin{defn}[Hartwig, Larsson and Silvestrov \cite{HLS}]
A Hom-Lie algebra is a triple $(V, [\cdot, \cdot], \alpha)$
consisting of
 a linear space $V$, bilinear map $[\cdot, \cdot]: V\times V \rightarrow V$ and
 a linear space homomorphism $\alpha: V \rightarrow V$
 satisfying
\begin{eqnarray} \label{skewsymHomLie}
 [x,y]=-[y,x] \quad && {\text{(skew-symmetry)}}
\\{}
\label{JacobyHomLie}
 \circlearrowleft_{x,y,z}{[\alpha(x),[y,z]]}=0
\quad && {\text{(Hom-Jacobi identity)}}
\end{eqnarray}
for all $x, y, z$ from $V$, where
$\circlearrowleft_{x,y,z}$ denotes summation over
the cyclic permutation on $x,y,z$.
\end{defn}
Using the skew-symmetry, one may write the
Hom-Jacobi identity in the form
(\ref{LeibnizalgidentAdDeriv}). Thus, if a
Hom-Leibniz algebra is skewsymmetric, then it is
a Hom-Lie algebra.

\begin{prop}
Every skew-symmetric bilinear map on a
$2$-dimensional linear space defines a Hom-Lie
algebra.
\end{prop}
\begin{proof}
The Hom-Jacobi identity \eqref{JacobyHomLie} is
satisfied for any triple $(x,x,y)$.
\end{proof}

If we add the condition that the linear map $\alpha$ is an algebra
homomorphism with respect to the bracket then we get restricted
Hom-Lie algebras which are a special case of Quasi-Hom-Lie
algebras. All these classes are a special case of the more general
Quasi-Lie algebras \cite{LS1,LS2}. In the following, we recall the
definition of Quasi-Lie algebras.

We denote by $\mathfrak{L}_{\K}(V)$  the set of
linear maps of the linear space $V$ over the
field $\K$.
\begin{defn} (Larsson, Silvestrov \cite{LS2}) \label{def:quasiLiealg}
A \emph{quasi-Lie algebra} is a tuple
$(V,[\cdot,\cdot]_V,\alpha,\beta,\omega,\theta)$ consisting of
\begin{itemize}
    \item $V$ is a linear space over $\mathbb{K}$,
    \item $[ \cdot,\cdot]_V:V\times V\to V$ is a bilinear
      map that is called a product or bracket in $V$;
    \item $\alpha,\beta:V\to V$ are linear maps,
    \item $\omega:D_\omega\to \mathfrak{L}_{\mathbb{K}}(V)$ and $\theta:D_\theta\to \mathfrak{L}_{\mathbb{K}}(V)$
      are maps with domains of definition $D_\omega, D_\theta\subseteq V\times
    V$,
\end{itemize}
such that the following conditions hold:
\begin{itemize}
      \item ($\omega$-symmetry) The product satisfies a generalized skew-symmetry condition
        $$[ x,y]_V=\omega(x,y)[ y,x]_V,
        \quad\text{ for all } (x,y)\in D_\omega ;$$
\item (quasi-Jacobi identity) The bracket satisfies a generalized Jacobi identity
    $$\circlearrowleft_{x,y,z}\big\{\,\theta(z,x)\big([\alpha(x),[ y,z]_V]_V+
    \beta [ x,[ y,z]_V ]_V\big)\big\}=0,$$
     for all $(z,x),(x,y),(y,z)\in D_\theta$.
\end{itemize}
\end{defn}
Note that $(\omega(x,y)\omega(y,x)-id)[ x,y]=0,$ if $(x,y), (y,x)
\in D_\omega$, which follows from the computation $[
x,y]=\omega(x,y)[ y,x] =\omega(x,y)\omega(y,x)[ x,y].$

The class of Quasi-Lie algebras incorporates as special cases
\emph{hom-Lie algebras} and more general \emph{quasi-hom-Lie
algebras (qhl-algebras)} which appear naturally in the algebraic
study of $\sigma$-derivations (see \cite{HLS}) and related
deformations of infinite-dimensional and finite-dimensional Lie
algebras. To get the class of qhl-algebras one specifies
$\theta=\omega$ and restricts attention to maps $\alpha$ and $\beta$
satisfying the twisting condition $[
\alpha(x),\alpha(y)]=\beta\circ\alpha [ x,y]$. Specifying this
further by taking $D_\omega =V \times V$, $\beta=id$ and
$\omega=-id$, one gets the class of Hom-Lie algebras including Lie
algebras when $\alpha = id$. The class of quasi-Lie algebras
contains also color Lie algebras and in particular Lie
superalgebras.

\begin{prop} \label{prop:functorHomassHomLie} {\bf (Functor Hom-Lie)}
To any Hom-associative algebra defined by the multiplication $\mu$
and a homomorphism $\alpha$ over a $\K$-linear space $V$, one may
associate a  Hom-Lie algebra defined for all $x,y \in V$ by the
bracket
$$
[ x,y ]=\mu (x,y)-\mu (y,x )
$$
\end{prop}
\begin{proof}
The bracket is obviously skewsymmetric and with a direct computation
we have
\begin{multline*}
[\alpha (x),[y,z]]+[\alpha (z),[x,y]]+
  [\alpha(y),[z,x]]= \\
  \mu(\alpha (x),\mu(y,z))-\mu(\alpha(x),\mu(z,y))-\mu(\mu(y,z),\alpha(x))+
\mu(\mu(z,y),\alpha(x))\\
+\mu(\alpha(z),\mu(x,y))-\mu(\alpha(z),\mu(y,x))
-\mu(\mu(x,y),\alpha(z))+\mu(\mu(y,x),\alpha(z))
\\
+\mu(\alpha (y),\mu(z,x))-\mu(\alpha
(y),\mu(x,z))-\mu(\mu(z,x),\alpha(y))+\mu(\mu(x,z),\alpha(y))=0
\end{multline*}

\end{proof}
\section{Hom-Lie-Admissible algebras}
The Lie admissible algebras were introduced by A.
A. Albert in 1948. Physicists attempted to
introduce this structure instead of Lie algebras.
For instance, the validity of Lie-Admissible
algebras for free particles is well known. These
algebras arise also in classical quantum
mechanics as a generalization of conventional
mechanics (see \cite{HYO}, \cite{MOS}). In this
section, we extend to Hom-algebras the classical
concept of Lie-Admissible algebras.
\begin{defn} \label{def:Lieadmis}
Let $\A$ be a Hom-algebra structure on $V$
defined by the multiplication $\mu$ and a
homomorphism $\alpha$.  Then $\A$ is said to be
Hom-Lie admissible algebra over $V$ if the
bracket defined for all $x,y \in V$ by
\begin{equation} \label{commutator}
[ x,y ]=\mu (x,y)-\mu (y,x )
\end{equation}
satisfies the Hom-Jacobi identity $
\circlearrowleft_{x,y,z}{[\alpha(x),[y,z]]}=0 $
for all $x,y,z\in V$.
\end{defn}
\begin{rem}
Since the commutator bracket \eqref{commutator}
is always skew-symmetric, it makes any Hom-Lie
admissible algebra into a Hom-Lie algebra.
\end{rem}
\begin{rem}
According to Proposition
\ref{prop:functorHomassHomLie},
any Hom-associative algebra is Hom-Lie admissible
with the same twisting map.
\end{rem}

It can be also checked that any Hom-Lie algebra with a
twisting $\alpha$ is Hom-Lie admissible with the
same twisting $\alpha$.

\begin{prop} \label{prop:functorHomLieHomLie}
Any Hom-Lie algebra $(V, [\cdot, \cdot], \alpha)$
is Hom-Lie admissible with the same twisting map
$\alpha$.
\end{prop}
\begin{proof}
For the commutator product $\langle x, y \rangle
= [x,y] - [y,x]$ one has
\begin{eqnarray*}
\langle x, y \rangle = [x,y] - [y,x] = -([y,x] -
[x,y])= -\langle y, x \rangle, \\
\circlearrowleft_{x,y,z} \langle \alpha(x),
\langle y , z \rangle\rangle =
\circlearrowleft_{x,y,z} ([\alpha(x),[y, z]] -
[\alpha(x),[z, y]]  \\ {}
-[[y,z],\alpha(x)] + [[z, y],\alpha(x)]) = \\
4 \circlearrowleft_{x,y,z} [\alpha(x),[y, z]] =
0.
\end{eqnarray*}
\end{proof}

\subsection{Lie-Santilli isotopies of associative and Lie algebras}
Already as early as in 1967, Santilli
\cite{Sant67}, considered the two-parametric
deformations (mutations) of the Lie commutator
bracket in an associative algebra,
$(A,B)=pAB-qBA$, where $p$ and $q$ are scalar
parameters and $A$ and $B$ are elements in the
associative algebra (typically algebra of
matrices or linear operators); and in 1978,
Santilli \cite{Sant78-1,Sant78-2,Sant78} extended
it to ``operator-deformations" of the Lie product
as follows $(A,B)=APB-BQA$, where $P$ and $Q$ are
fixed elements in the underlying associative
algebra. The motivation for introducing
non-associative generalizations of the bracket
multiplication came from attempts to resolve
certain limitations of conventional formalism of
classical and quantum mechanics. Subsequently, in
numerous works including articles and books by
Santilli and other authors the evolution
equations based on such deformed brackets,
physical applications and physical consequences
of introducing such generalized models has been
investigated. The deformations of the commutator
bracket multiplication introduced by Santilli in
investigations of foundations of classical and
quantum mechanics and hadronic physics, have
reappeared later in many incarnations both in
Mathematics and Physics. We refer the reader to
\cite{Falcon2001,Sant78-1,Sant78-2,Sant78,Sant79,SantBook78,Sant97,
Sant99,Sant2006,Tsagas93} for further discussion
and references on the models based on
introduction of such non-associative deformed
commutator bracket multiplications instead of
commutator (Lie) bracket multiplication, their
bi-module type generalizations (genoalgebras) as
well as for a review of relation with other
appearances of so deformed commutator brackets in
physics and mathematics, for example in contexts
related to quantum algebras, quantum groups, Lie
algebras and superalgebras, Jordan and other
classes of algebras.

The relation between Hom-Lie admissible algebras
and non-associative algebras with Santilli's
deformed bracket products is certainly an
interesting direction for further investigation.
Here we would like to highlight some
interrelations and differences between considered
algebraic objects.

The Santilli's products with arbitrary $P$ and
$Q$ are not anti-symmetric in general except when
$P$ and $Q$ are specially interrelated within the
underlying algebra, for example when $P=Q$ that
is $(A,B)=AQB-BQA$. It can also be checked that
the Lie-Santilli bracket product $(A,B)=AQB-BQA$
satisfies the Lie algebras Jacobi identity, which
is not the case in general for $P\neq Q$. Thus
the associative algebras with the modified
product $A\times_Q B = AQB$ are Lie admissible
algebras. Since $$<A,B>=(A,B)-(B,A)=
A(P+Q)B-B(P+Q)A,$$ the Santilli's deformed
commutator bracket product $(A,B)=APB-BQA$
defines a Lie admissible algebra as well.

\begin{rem}
The Hom-Lie admissible algebras are Lie
admissible algebras typically when the twisting
mapping $\alpha$ built into their definition is
the scalar multiple of the identity mapping.
It could thus be of interest to describe, for a given twisting
$\alpha$, as sharp as possible conditions on $P$ and $Q$ in order for
the general Santilli's bracket product
$(A,B)=APB-BQA$ to define a Hom-Lie algebra or more
general Hom-Leibniz algebra with twisting $\alpha$.
\end{rem}
\begin{rem}
It is possible also that Hom-Lie
admissible algebras can be studied
using Santilli's formalism of genoalgebras, the abstract bi-module type
extension of the deformed commutator brackets $(A,B)=APB-BQA$, but after this formalism
is appropriately modified to suit the Hom-algebras context.
\end{rem}

Santilli  has considered also so called isotopies
of associative and Lie algebras. Algebraic
problem can be formulated as follows. Any
associative algebra is Lie admissible since the
commutator bracket on any associative algebra
satisfies axioms of a Lie algebra. How can
associative products in associative algebras be
modified to yield as general as possible Lie
admissible algebras? Can any Lie admissible
algebra be obtained by such modifications from
some associative algebra?  Such modifications of
associative and corresponding Lie algebras where
called isotopies of associative and Lie algebras.
In
\cite{Sant78-1,Sant78-2,Sant78,Sant79,Sant97,Sant99,Sant2006},
several general isotopies of associative products
and associated Lie products have been identified.
The most general of all presented there isotopies
of a product for elements in an associative
algebra $A$ over a field $\K$ is given by $x*y =
awxwTwyw$, where $a\in \K, w\in A, w^2=w \neq 0$,
and $T$ is some extra element. The product $*$ is
associative if $w^2=w$ and $T\in A$. Santilli
allows $T$ to be some extra element outside $A$.
Then a special care is needed on algebraic side
in order to make involved objects and maps to be
properly defined. If $a(w)(x)(w)T(w)(y)(w)$ is
not identified with some element of $A$, then the
new product $*$ is taking values in some
generally non-associative algebra $A_T$
generated, as a linear space over $\K$, by
elements $x\in A$ and formal expressions of the
form $x_1Tx_2 Tx_3\dots x_{n-1} T x_{n} \in ATATA
\dots ATA$ for $x_1, \dots , x_{n}\in A$ for
integers $n\geq 2$, whatever these expressions
mean. With a $\K$-bilinear product on $A_T$
defined for $x, x_1, \dots , x_{n}, y_1, \dots ,
y_{m}\in A$ by
\begin{eqnarray*}
&& (x_1Tx_2Tx_3 \dots x_{n-1} T x_{n})
(y_1Ty_2Ty_3 \dots y_{n-1} T y_{m}) = \\
&& x_1Tx_2Tx_3 \dots x_{n-1} T (x_{n}
y_1)Ty_2Ty_3 \dots y_{m-1} T y_{m}, \\
&& x (x_1Tx_2Tx_3 \dots x_{n-1} T x_{n}) = (xx_1)Tx_2Tx_3 \dots x_{n-1} T x_{n}, \\
&& ( x_1Tx_2Tx_3 \dots x_{n-1} T x_{n})x  = x_1Tx_2Tx_3 \dots
x_{n-1} T (x_{n}x),
\end{eqnarray*}
in particular, $u(xTy)= (ux)T(y)$, $(xTy)v = (x)T(yv)$, $u(xTy)v =
(ux)T(yv)$ hold for $u,x,y,v \in A$. Then the expression $x*y =
a(w)(x)(w)T(w)(y)(w)$ yields again the element from $A_T$ for $x, y
\in A_T$, and we get the product $*$ on the algebra $A_T$. If now
$w^2=w$, then $*$ satisfies the associativity condition $x*(y*z) =
(x*y)*z$ on $A_T$. Indeed, for $x=x_1Tx_2Tx_3 \dots x_{n-1} T
x_{n}$, $y=y_1Ty_2Ty_3 \dots y_{n-1} T y_{m}$ and $z=z_1Tz_2Tz_3
\dots z_{n-1} T z_{k}$, using $w^2=w$, one gets
\begin{eqnarray*}
x*(y*z) &=& a(w)( x_1Tx_2Tx_3 \dots x_{n-1} T x_{n})(w) T (w)
\\
&&(a(w)( y_1Ty_2Ty_3 \dots y_{n-1} T y_{m})(w)T(w) \\
&& ( z_1Tz_2Tz_3 \dots z_{n-1} T z_{k})(w))(w)
\\
&=& a^2(w)( x_1Tx_2Tx_3 \dots x_{n-1} T x_{n})(w) T
\\
&&(w) ( y_1 Ty_2Ty_3 \dots y_{n-1} T y_{m})(w)T(w) \\
&& ( z_1Tz_2Tz_3 \dots z_{n-1} T z_{k})(w)
\\
&=& a(w)(a(w)( x_1Tx_2Tx_3 \dots x_{n-1} T x_{n})(w) T (w)
\\
&&( y_1 Ty_2Ty_3 \dots y_{n-1} T y_{m})(w))(w)T(w) \\
&& ( z_1Tz_2Tz_3 \dots z_{n-1} T z_{k})(w)
\\
&=&(x*y)*z.
\end{eqnarray*}
One may show that this algebra carries a
structure of Hom-associative algebra  in the
following way. Let $\alpha (s)=wsw$ for all $s\in
A_T$. Then
\begin{eqnarray*}
&& \alpha(x)*(y*z) = aw (wxw) w T w (awywTwzw)w  \\
&&= awxwTw(awywTwzw)w  \\
&&= x*(y*z) = (x*y)*z  \\
&&= aw(awxwTwyw)wTwzw \\
&&= aw(awxwTwyw)wT w (wzw) w \\
&&= (x*y)*\alpha(z)
\end{eqnarray*}
The map $\alpha: s\mapsto
wsw$ is a linear map on $A_T$, and $(A_T,*)$ is at the same time
Hom-associative with non-trivial twisting map $\alpha$. However,
since the product $*$ on $A_T$ is associative, the algebra $(A_T,*)$
is Lie admissible, or in other words a Hom-Lie admissible with the
trivial twisting map $id_{A_T}$.

If $w^2 \neq w$, then the above reduction of Hom-associativity with
the twisting map $\alpha:s\mapsto wsw$ to associativity is not
working. Moreover, on many associative algebras $A$, there are
linear maps $\alpha$, which can not be represented on the form
$x\mapsto uxv$ at all. This explains why the classes of
Hom-associative and Hom-Lie (admissible) algebras are different from
Lie-Santilli algebras and isotopies of associative and Lie algebras.
Whereas the relation between Hom-associative and Hom-Lie
(admissible) algebras resembles that between associative and Lie
algebras, the Hom-associative algebras are mostly not associative,
and Hom-Lie (admissible) algebras are mostly not Lie (admissible)
algebras. Only with special choices of $\alpha$ like above one
recovers the associative and Lie algebra properties holding also in
the case of Santilli's isotopies of associative and Lie algebras.

\begin{rem}
We assumed that $w\in A$ for simplicity of exposition. All
conclusions above hold however even if $w$ as $T$ is an extra
element possibly outside $A$, but of cause with $A_T$ replaced by a
properly defined algebra $A_{T,w}$ built from elements $x\in A$ and the formal
product expressions $xw$, $wx$, $xT$, $Tx$, $Tw$, $wT$, with
multiplication obeying the same rules as in $A_{T}$, allowing to
group the elements in the products by brackets in the appropriate
way, as if $w$ would belonged to $A$.
\end{rem}

The other fundamental algebraic issue tackled by Santilli is imbedding
of the scalar field into the algebras over this filed. If an algebra
$A$ over the field $\K$ with a unit $1_\K$ has a unit $1_A$, then
there is a canonical imbedding of the field into the algebra given
by $i_A:c \mapsto c 1_A$ for $c\in \K$.
%
Also one has  $1_A x=x1_A=x$. If the multiplication in $A_T$ is
defined as $x*y=xTy$ (corresponding to $a=1_\K$ and $w=1_A$) then one still would like to have $1_{A_T}
*x=x*1_{A_T}=x$, which can be written as
$1_{A_T} T x=x T 1_{A_T}=x$. Thus, if $T$ has
left and right inverse $T^{-1}$, then
$1_{A_T}=T^{-1}$. Also the canonical imbedding of
the field into the new algebra yields
$i_{A_T}(1_\K)= 1_\K 1_{A_T} = 1_\K T^{-1}$ and
more generally $i_{A_T}(c)=c 1_{A_T} = c T^{-1}$
for $c\in \K$. These elements form a field inside
$A_T$ with the unit $\hat{I}=1_\K T^{-1}$. This
field $\hat{\K}$ is called isofield. Santilli
have noticed that dependence on $T$ of the new
unit in the isofield, caused by the changed
product in the algebra, is not just some
complicating curiosity, but advantageous
phenomena that opens new vast fundamental
opportunities in physics, differential geometry,
tensor calculus and beyond. This is because,
while the unit $1_\K$ in the scalar field $\K$ is
fixed, $T$ and thus the unit $\hat{I}$ in the
isofield $\hat{\K}$ can be chosen to depend on
the non-linear functionals or expressions in some
other parameters, functions, their derivatives,
integrals, etc., in physics having interpretation
as time, position, speed, momentum, acceleration,
mass, energy, etc. This dependence may be well
highly non-linear. Santilli have made an effort
in systematic analysis of how the new algebra
structures and introduction of isofield
$\hat{\K}$ and parameter dependent isounits
effect the equations of motion, time evolution
and other basics of Hamiltonian and quantum
mechanics. These first steps open a huge field
for further research in many directions of
interest both in physics and mathematics.

%

\section{Classification of Hom-Lie admissible
algebras} In the following, we explore some other general classes of
Hom-Lie admissible algebras, $G$-Hom-associative algebras, extending the class of Hom-associative algebras.
Is is convenient first to introduce the
following notation.
\begin{defn}
By $\alpha$-associator of $\mu$ we call a trilinear map
$a_{\alpha,\mu}$ over $V$ associated to a product $\mu$ and a
homomorphism $\alpha$ defined  by
$$
a_{\alpha,\mu}(x_1,x_2,x_3)=\mu (\mu (x_1,x_2),\alpha
(x_{3}))-\mu(\alpha(x_{1}),\mu (x_{2},x_3)).
$$
\end{defn}
\begin{defn}
Let $G$ be a subgroup of the permutations group $\mathcal{S}_3$, a
Hom-algebra on $V$ defined by the multiplication $\mu$ and a
homomorphism $\alpha$ is said $G$-Hom-associative if
\begin{equation}\label{admi}
\sum_{\sigma\in G} {(-1)^{\varepsilon ({\sigma})}
(\mu (\mu (x_{\sigma (1)},x_{\sigma (2)}),\alpha
(x_{\sigma (3)}))}-\mu(\alpha(x_{\sigma (1)}),\mu
(x_{\sigma (2)},x_{\sigma (3)}))=0
\end{equation}
where $x_i$ are in $V$ and $(-1)^{\varepsilon ({\sigma})}$ is the
signature of the permutation $\sigma$.
\end{defn}
The condition (\ref{admi}) may be written
\begin{equation}
\sum_{\sigma\in G}{(-1)^{\varepsilon ({\sigma})}a_{\alpha,\mu}\circ
\sigma}=0
\end{equation}
where
$\sigma(x_1,x_2,x_3)=(x_{\sigma(1)},x_{\sigma(2)},x_{\sigma(3)})$.

\begin{rem}
If $\mu$ is the multiplication of a Hom-Lie
admissible Lie algebra, then the condition
\eqref{admi} is equivalent to the property that
the bracket defined by
$$
[ x,y ]=\mu (x,y)-\mu (y,x )
$$
satisfies the Hom-Jacobi identity, or
equivalently
\begin{equation}
\sum_{\sigma\in \mathcal{S}_3}
{(-1)^{\varepsilon({\sigma})} (\mu (\mu
(x_{\sigma (1)},x_{\sigma (2)}),\alpha (x_{\sigma
(3)}))}-\mu(\alpha(x_{\sigma (1)}),\mu (x_{\sigma
(2)},x_{\sigma (3)})))=0
\end{equation}
which may be written as
\begin{equation}
\sum_{\sigma\in \mathcal{S}_3}{(-1)^{\varepsilon
({\sigma})}a_{\alpha,\mu}\circ \sigma}=0.
\end{equation}
\end{rem}

\begin{prop}
Let $G$ be a subgroup of the permutations group
$\mathcal{S}_3$. Then any $G$-Hom-associative
algebra is a Hom-Lie admissible algebra.
\end{prop}
\begin{proof}The skewsymmetry follows straightaway from the
definition.

We have a subgroup $G$ in $\mathcal{S}_3$. Take the set of conjugacy
class $\{g G\}_{g\in I}$  where $I\subseteq G$, and for any
$\sigma_1, \sigma_2\in I,\sigma_1 \neq \sigma_2 \Rightarrow \sigma_1
G\bigcap \sigma_1 G =\emptyset$. Then

$$\sum_{\sigma\in \mathcal{S}_3}{(-1)^{\varepsilon ({\sigma})}a_{\alpha,\mu}\circ
\sigma}=\sum_{\sigma_1\in I}{\sum_{\sigma_2\in \sigma_1
G}{(-1)^{\varepsilon ({\sigma_2})}a_{\alpha,\mu}\circ \sigma_2}}=0$$
\end{proof}

The $G$-associative algebra in classical case was studied in
\cite{GR04}. The result may be extended to Hom-structures in the
following way.

 The subgroups of $\mathcal{S}_3$ are
$$G_1=\{Id\}, ~G_2=\{Id,\tau_{1 2}\},~G_3=\{Id,\tau_{2
3}\},$$ $$~G_4=\{Id,\tau_{1 3}\},~G_5=A_3 ,~G_6=\mathcal{S}_3,$$
where $A_3$ is the alternating group and where $\tau_{ij}$ is the
transposition between $i$ and $j$.

We obtain the following type of Hom-Lie
admissible algebras.
\begin{itemize}
\item The  $G_1$-Hom-associative algebras  are the Hom-associative
algebras defined above.

\item The  $G_2$-Hom-associative algebras satisfy the condition
\begin{equation}\nonumber
\mu(\alpha(x),\mu (y,z))-\mu(\alpha(y),\mu (x,z))=\mu (\mu
(x,y),\alpha (z))-\mu (\mu (y,x),\alpha (z))
\end{equation}
When $\alpha$ is the identity the algebra is called Vinberg algebra
or left symmetric algebra.
\item The  $G_3$-Hom-associative algebras satisfy the condition
\begin{equation}\nonumber
\mu(\alpha(x),\mu (y,z))-\mu(\alpha(x),\mu (z,y))=\mu (\mu
(x,y),\alpha (z))-\mu (\mu (x,z),\alpha (y))
\end{equation}
When $\alpha$ is the identity the algebra is called pre-Lie algebra
or right symmetric algebra.
\item The  $G_4$-Hom-associative algebras satisfy the condition
\begin{equation}\nonumber
\mu(\alpha(x),\mu (y,z))-\mu(\alpha(z),\mu (y,x))=\mu (\mu
(x,y),\alpha (z))-\mu (\mu (z,y),\alpha (x))
\end{equation}
\item The  $G_5$-Hom-associative algebras satisfy the condition
\begin{eqnarray}\nonumber
\nonumber \mu(\alpha(x),\mu (y,z))+\mu(\alpha(y),\mu
(z,x)+\mu(\alpha(z),\mu
(x,y))= \\
\nonumber \mu (\mu (x,y),\alpha (z))+\mu (\mu (y,z),\alpha (x))+\mu
(\mu (z,x),\alpha (y))
\end{eqnarray}\nonumber
If the product $\mu$ is skewsymmetric then the previous condition is
exactly the Hom-Jacobi identity.
\item The  $G_6$-Hom-associative algebras are the Hom-Lie admissible
algebras.
\end{itemize}
Then, one may set the following generalization of Vinberg and
pre-Lie algebras.
\begin{defn}

A Hom-Vinberg algebra is a triple $(V, \mu, \alpha)$ consisting of a
linear space $V$, a bilinear map $\mu: V\times V \rightarrow V$ and
a homomorphism $\alpha$ satisfying
\begin{equation}
\mu(\alpha(x),\mu (y,z))-\mu(\alpha(y),\mu (x,z))=\mu (\mu
(x,y),\alpha (z))-\mu (\mu (y,x),\alpha (z))
\end{equation}
\end{defn}

\begin{defn}
A Hom-pre-Lie  algebra is a triple $(V, \mu, \alpha)$ consisting of
a linear space $V$, a bilinear map $\mu: V\times V \rightarrow V$
and a homomorphism $\alpha$ satisfying
\begin{equation}
\mu(\alpha(x),\mu (y,z))-\mu(\alpha(x),\mu (z,y))=\mu (\mu
(x,y),\alpha (z))-\mu (\mu (x,z),\alpha (y))
\end{equation}
\end{defn}

\begin{rem}
A Hom-pre-Lie algebra is the opposite algebra of a Hom-Vinberg
algebra.
\end{rem}
\section{Flexible Hom-Lie admissible algebras}
  The study of flexible Lie admissible algebras was initiated by Albert
  \cite{Albert48} and investigated by number of authors Myung,
  Okubo, Laufer, Tomber and Santilli,  see \cite{myung78} and
  \cite{Sant79}. The aim of this section is
  to extend the notions and results about flexible
  Lie admissible algebras to Hom-structures.
 \begin{defn}
 A Hom-algebra $\A =(V, \mu, \alpha)$ is called
 flexible if    for any $x,y$ in $V$
\begin{equation}\label{flexible}
      \mu (\mu (x,y), \alpha (x))=   \mu (\alpha (x),\mu (y,x)))
\end{equation}

 \end{defn}
 \begin{rem}
 Using the $\alpha$-associator
 $$  a_{\mu , \alpha}(x,y,z)=
 \mu(\mu(x,y),\alpha(z))-\mu(\alpha(x),\mu(y,z)),
 $$
  the condition (\ref{flexible}) may be written as
 \begin{equation}\label{flexible2}
      a_{\mu , \alpha}(x,y,x)=0
\end{equation}
\end{rem}

The $\alpha$-associator $a_{\mu , \alpha}$ is a
useful tool in the study of Hom-Lie and Hom-Lie
admissible algebras.

\begin{lem}\label{flexi23}
Let $\A =(V, \mu, \alpha)$ be a Hom-Lie
admissible algebra. The following assertions are
equivalent
\begin{enumerate}
\item $\A$ is flexible.
\item For any $x,y$ in $V$, $a_{\mu , \alpha}(x,y,x)=0$.
\item For any $x,y,z$ in $V$,
$a_{\mu , \alpha}(x,y,z)=-a_{\mu
,\alpha}(z,y,x).$
\end{enumerate}
\end{lem}
\begin{proof}
The equivalence of the first two assertions
follows from the definition, and of last two
assertions since $ a_{\mu , \alpha}(x-z,y,x-z)=0
$ holds by definition and is equivalent to $
a_{\mu , \alpha}(x,y,z)+ a_{\mu,\alpha}(z,y,x)=0$
by linearity.
\end{proof}
\begin{cor}
Any Hom-associative algebra is flexible.
\end{cor}
In the following, we aim to characterize  the
flexible Hom-Lie admissible algebras.

Let $\A =(V, \mu, \alpha)$ be a Hom-algebra  and
let $[x,y]= \mu(x,y)-\mu(y,x) $ be its
commutator. We introduce the notation
$$S(x,y,z):= a_{\mu , \alpha}(x,y,z)+
a_{\mu , \alpha}(y,z,x)+a_{\mu , \alpha}(z,x,y).
$$
Then we have the following properties.
\begin{lem}
$$ S(x,y,z)=[\mu(x,y),\alpha(z)]+
[\mu(y,z),\alpha(x)]+[\mu(z,x),\alpha(y)].
$$
\end{lem}
\begin{proof} The assertion follows by expanding
commutators on the right hand side:
$$[\mu(x,y),\alpha(z)]+ [\mu(y,z),\alpha(x)]+[\mu(z,x),\alpha(y)]=$$
$$\mu(\mu(x,y),\alpha(z))-\mu(\alpha(z),\mu(x,y))+\mu(\mu(y,z),\alpha(x))-\mu(\alpha(x),\mu(y,z))+$$
$$\mu(\mu(z,x),\alpha(y))-\mu(\alpha(y),\mu(z,x))
=S(x,y,z).$$
\end{proof}
\begin{prop}
A Hom-algebra $\A$ is Hom-Lie admissible if and
only if it satisfies
$$ S(x,y,z)=S(x,z,y)$$
for any  $x,y,z \in V.$
\end{prop}
\begin{proof} The assertion follows from
$$S(x,z,y)-S(x,y,z)=$$
$$[\mu(x,z),\alpha(y)]+ [\mu(z,y),\alpha(x)]+
[\mu(y,x),\alpha(z)]$$
$$-[\mu(x,y),\alpha(z)]-
[\mu(y,z),\alpha(x)]-[\mu(z,x),\alpha(y)]=$$
$$\circlearrowleft_{x,y,z}[\alpha(x),[y,z]].$$
\end{proof}

Let $\A =(V,\mu , \alpha)$ be a Hom-algebra where
 $\mu$ is the multiplication and  $\alpha$
 a homomorphism. We denote by
 $\A ^+$ the Hom-algebra over $V$ with a multiplication $x \bullet
 y=\frac{1}{2}(\mu (x,y) +\mu(y,x))$.
 We also denote by   $\A ^-$ the Hom-algebra
 over $V$ where the
 multiplication is given by the commutator
 $[x,y]=\mu (x,y) -\mu(y,x)$.
\begin{prop} A Hom-algebra $\A =(V,\mu , \alpha)$
is flexible if and only if
$$ \label{equa2}
[\alpha(x),y\bullet z]=
[x,y]\bullet\alpha(z)+\alpha(y)\bullet[x,z].
$$
\end{prop}
\begin{proof}
Let $\A$ be a flexible Hom-algebra. By Lemma
\ref{flexi23}, this is equivalent to $a_{\mu
,\alpha}(x,y,z)+ a_{\mu ,\alpha}(z,y,x)=0$ for
any $x,y,z$ in $V$, where $a_{\mu , \alpha}$ is
the $\alpha$-associator of $\A$.
 This implies
\begin{eqnarray}\label{equa} a_{\mu , \alpha}(x,y,z)+a_{\mu , \alpha}(z,y,x)+a_{\mu ,
 \alpha}(x,z,y)+a_{\mu , \alpha}(y,z,x)\\ \nonumber -a_{\mu , \alpha}(y,x,z)-a_{\mu ,
 \alpha}(z,x,y)=0.
 \end{eqnarray}
By expansion, the previous relation is equivalent to
$$
  [\alpha(x),y\bullet z]=[x,y]\bullet\alpha(z)+\alpha(y)\bullet[x,z]
$$
Conversely, assume we have the condition in
Proposition \ref{equa2}. By setting $x=z$ in
\eqref{equa}, one gets $a_{\mu ,
\alpha}(x,y,x)=0$. Therefore $\A$ is flexible.
\end{proof}
\section{Algebraic varieties of Hom-structures}
Let $V$ be a $n$-dimensional $\K$-linear space and
$\{e_1,\cdots,e_n\}$ be a basis of $V$. A  Hom-algebra structure on
$V$ with product $\mu$ is determined by $n^3$ structure constants
$C_{i j}^k$, where $\mu (e_i,e_j)=\sum_{k=0}^{n}{C_{i j}^k e_k}$ and
homomorphism $\alpha$ which is given by $n^2$ structure constants
$a_{i j}$, where $\alpha (e_i)=\sum_{j=0}^{n}{a_{i j} e_j}$.

If we require this algebra structure to be
Hom-associative, then this limits the set of
structure constants $(C_{i j}^k , a_{i j})$ to a
subvariety of $\K ^{n^3+n^2}$ defined by the
following polynomial equations:
$$
\sum_{l,m=1}^{n}{a_{i l} C_{j k}^{m} C_{l m}^{s}-
a_{k m}C_{i j}^{l} C_{l m}^{s}}=0, \quad
i,j,k,s=1,\cdots,n.
$$
The algebraic variety of $n$-dimensional Hom-associative algebras is
denoted by $HomAss_n$. Note that the equations are given by cubic
polynomials.

If we consider the $n$-dimensional unital
Hom-associative algebras, then we obtain a
subvariety which we denote by $HomAlg_n$ and
determined by the following polynomial equations:
\begin{eqnarray}
\begin{array}{l}
    \sum_{l,m=1}^{n}{a_{i l} C_{j k}^{m} C_{l m}^{s}- a_{k m}C_{i j}^{l}
C_{l m}^{s}}=0 \\[0.3cm]
 C_{i 1}^k=C_{1 i}^k=\left\{
 \begin{array}{l} 1,\ \text{if} \ i=k \\
 0,\ \text{if} \ i\neq k
 \end{array} \right.
\quad \quad \text{ for } i,j,k,s=1,\cdots,n.
\end{array}
\end{eqnarray}

 If we require that this algebra structure to be Hom-Lie,
then the structure constants $\{(C_{i j}^k)_{i<j} , (a_{i j})\}$
determine a subvariety of $\K ^{n^2(n+1)/2}$ defined by the
following system of polynomial equations
$$
\sum_{l,m=1}^{n}{a_{i l} C_{j k}^{m} C_{l m}^{s}+ a_{j l}C_{k i}^{m}
C_{l m}^{s}+ a_{k l} C_{i j}^{m} C_{l m}^{s}}=0
$$
where $C_{i j}^{k}=-C_{j i}^{k}.$

The variety of $n$-dimensional Hom-Lie algebras is denoted by
$HomLie_n$.

One can do the same for every Hom-Lie admissible
structure. If the algebras are parameterized by
their structure constants and the structure
constants of the homomorphism, the set of
$n$-dimensional algebras of given structure
carries a structure of algebraic varieties in
$\K^{n^3+n^2}$ defined by a system of cubic
polynomials.

A point in such an algebraic variety (e.g. , $HomAss_n$, $HomLie_n$)
represents an $n$-dimensional algebra (e.g. Hom-associative, or
Hom-Lie), along with a particular choice of basis. A change of basis
may give rise to a (possibly) different point of the algebraic
variety. The group $GL(n,\K)$ acts on the algebraic varieties of
Hom-structures by the so-called "transport of structure" action
defined as follows :

Let $\A$ be an algebra (e.g. Hom-associative, or Hom-Lie) defined by
multiplication $\mu$ and homomorphism $\alpha$. Given  $f\in
GL(n,\K)$, the action $f\cdot \A$ transports the structures by
\begin{eqnarray}
f\cdot \mu (x,y)=f^{-1}\mu (f(x),f(y))\\
f\cdot \alpha (x)=f^{-1}\alpha (f(x))
\end{eqnarray}
The \textit{orbit} of an algebra (e.g.
Hom-associative, or Hom-Lie) $\A$ is given by
$$\vartheta \left( \A\right) =\left\{ \A^{\prime
}=f\cdot \A,\quad f\in GL_n\left( \K\right)
\right\} .$$ The orbits are in 1-1-correspondence
with the isomorphism classes of $n$-dimensional
algebras.

\section{Three-dimensional Hom-Lie algebras}
In this section, we provide a characterization of all Hom-Lie
algebras of $sl(2)$ type and construct all the $3$-dimensional
Hom-Lie algebras for a class of homomorphisms.

 By a direct
calculation, we obtain the following class of Hom-Lie algebras of
$sl(2)$ type.
\begin{prop}
Let $V$ be a three dimensional $\K$-linear space and let $(X_1 ,X_2,
X_3)$ be its basis. Any Hom-Lie algebra with the following brackets
$$[X_1,X_2]=2 X_2, \ [X_1,X_3]=-2 X_3, \ [X_2,X_3]=X_1
$$ is given by homomorphism $\alpha$ defined, with respect to the previous basis, by a matrix of the form
 $\left(
  \begin{array}{ccc}
  a & d& c \\
  2c & b & f \\
    2d & e & b\\
  \end{array}
\right) $ where $a,c,d,e,f$ are any parameters in $\K$

\end{prop}
This may be viewed as a quasi-deformation of
$sl(2)$ in the class of Hom-Lie algebras. In case
where the matrix is the identity matrix one gets
the classical Lie algebras $sl(2)$.

\begin{rem}
It would be interesting to compare this list with
the quasi-Hom-Lie and Hom-Lie algebras of
quasi-deformations of $sl(2)$ constructed in
\cite{LS3} using a quasi-Hom-Lie structure on
twisted vector fields.
\end{rem}

We provide in the following a list of all Hom-Lie algebras
associated to the homomorphism  given  with respect to the basis
${e_1,e_2,e_3}$ of the $3$-dimensional $V$ by the following matrix
 $\left(
  \begin{array}{ccc}
  a & 0 & 0 \\
    0 & b & 0 \\
    0 & 0 & b\\
  \end{array}
\right) $.

In the following list, $C_{i j}^{k}$ for $i,j,k=1,2,3$ are
parameters in $\K$.

$\begin{array}[t]{l}
 [e_1,e_2 ]_1=C_{1 2}^{1}\ e_1+ C_{1 2}^{3}\ e _3\\{}
 [e_1,e_3 ]_1=C_{1 3}^{2}\ e_2\\{}
 [e_2,e_3 ]_1=C_{2 3}^{1} \ e_1+\frac{b \ C_{1 2}^{1}}{a}\ e_3\\
\end{array}$
\smallskip

$\begin{array}[t]{l}
 [e_1,e_2 ]_2 =\frac{a\ C_{2 3}^{3}}{b}\ e_1 - C_{1 3}^{3}\ e_2 + C_{1 2}^{3}\ e_3 \\{}
 [e_1,e_3 ]_2 =-\frac{a\ C_{2 3}^{2}}{b}\ e_1 + C_{1 3}^{2}\ e_2+C_{1 3}^{3}\ e_3 \\{}
 [e_2,e_3 ]_2 =C_{2 3}^{1}\ e_1+ C_{2 3}^{2}\ e _2+C_{2 3}^{3}\ e_3 \\
\end{array}$
\smallskip

  $\begin{array}[t]{l}
 [e_1,e_2 ]_3 =0 \\{}
 [e_1,e_3 ]_3 = C_{1 3}^{1}\ e_1 + C_{1 3}^{2}\ e_2 \\{}
 [e_2,e_3 ]_3 =C_{2 3}^{1}\ e_1\\
\end{array}$
\smallskip

 $\begin{array}[t]{l}
 [e_1,e_2 ]_4 = C_{1 2}^{3}\ e_3 \\{}
 [e_1,e_3 ]_4 = C_{1 3}^{2}\ e_2+C_{1 3}^{3}\ e_3 \\{}
 [e_2,e_3 ]_4 =0\\
\end{array}$
\smallskip

 $\begin{array}[t]{l}
 [e_1,e_2 ]_5 =C_{1 2}^{3}\ e_3 \\{}
 [e_1,e_3 ]_5 = C_{1 3}^{2}\ e_2 \\{}
 [e_2,e_3 ]_5 =C_{2 3}^{1}\ e_1\\
\end{array}$
\smallskip

 $\begin{array}[t]{l}
 [e_1,e_2 ]_6 = C_{1 2}^{1}\ e_1  + C_{1 2}^{3}\ e_3 \\{}
 [e_1,e_3 ]_6 =0 \\{}
 [e_2,e_3 ]_6 =C_{2 3}^{1}\ e_1\\
\end{array}$
 \smallskip

$\begin{array}[t]{l}
 [e_1,e_2 ]_7 = C_{1 2}^{1}\ e_1 + C_{1 2}^{2}\ e_2 + C_{1 2}^{3}\ e_3 \\{}
 [e_1,e_3 ]_7 =0\\{}
 [e_2,e_3 ]_7 =C_{2 3}^{1}\ e_1+C_{2 3}^{3}\ e_3\\
\end{array}$
\smallskip

 $\begin{array}[t]{l}
 [e_1,e_2 ]_8 = - C_{1 3}^{3}\ e_2 + C_{1 2}^{3}\ e_3 \\{}
 [e_1,e_3 ]_8 =C_{1 3}^{2}\ e_2+C_{1 3}^{3}\ e_3\\{}
 [e_2,e_3 ]_8 =C_{2 3}^{1}\ e_1\\
\end{array}$
\smallskip

$\begin{array}[t]{l}
 [e_1,e_2 ]_{9} = - C_{1 3}^{3}\ e_2 + C_{1 2}^{3}\ e_3 \\{}
 [e_1,e_3 ]_{9} =-\frac{a\ C_{2 3}^{2}}{b}\ e_1+C_{1 3}^{2}\ e_2+C_{1 3}^{3}\ e_3\\{}
 [e_2,e_3 ]_{9} =C_{2 3}^{1}\ e_1+C_{2 3}^{2}\ e _2\\
\end{array}$
\smallskip

$\begin{array}[t]{l}
 [e_1,e_2 ]_{10} =0 \\{}
 [e_1,e_3 ]_{10} = C_{1 3}^{1}\ e_1 +C_{1 3}^{2}\ e_2+C_{1 3}^{3}\ e_3 \\{}
 [e_2,e_3 ]_{10} =C_{2 3}^{1}\ e_1+C_{2 3}^{2}\ e _2\\
\end{array}$
\smallskip

$\begin{array}[t]{l}
 [e_1,e_2 ]_{11} =  C_{1 2}^{2}\ e_2 + C_{1 2}^{3}\ e_3 \\{}
 [e_1,e_3 ]_{11} =C_{1 3}^{2}\ e_2+C_{1 3}^{3}\ e_3\\{}
 [e_2,e_3 ]_{11} =0\\
\end{array}$
\smallskip

 $\begin{array}[t]{l}
 [e_1,e_2 ]_{12} =-\frac{C_{1 3}^{1}C_{2 3}^{3}}{C_{2 3}^{2}} \ e_1 +  C_{1 2}^{2}\ e_2 -\frac{C_{1 3}^{3}C_{2 3}^{3}}{C_{2 3}^{2}}\ e_3 \\{}
 [e_1,e_3 ]_{12} =C_{1 3}^{1}\ e_1 -\frac{C_{1 2}^{2}C_{2 3}^{2}}{C_{2 3}^{3}} \  e_2+C_{1 3}^{3}\ e_3\\{}
 [e_2,e_3 ]_{12} =C_{2 3}^{1}\ e_1+C_{2 3}^{2}\ e _2+C_{2 3}^{3}\ e_3\\
\end{array}$
\smallskip

$\begin{array}[t]{l}
 [e_1,e_2 ]_{13} =   C_{1 2}^{3}\ e_3 \\{}
 [e_1,e_3 ]_{13} =-\frac{a\ C_{2 3}^{2}}{b}\ e_1 +C_{1 3}^{2}\ e_2\\{}
 [e_2,e_3 ]_{13} =C_{2 3}^{1}\ e_1+C_{2 3}^{2}\ e _2\\
\end{array}$
\smallskip

$\begin{array}[t]{l}
 [e_1,e_2 ]_{14} =C_{1 2}^{1}\ e_1 -C_{1 3}^{3}\ e_2 + \frac{C_{1 2}^{1}C_{1 3}^{3}}{C_{1 3}^{1}}\ e_3 \\{}
 [e_1,e_3 ]_{14} = C_{1 3}^{1}\ e_1 - \frac{C_{1 3}^{1}C_{1 3}^{3}}{C_{1 2}^{1}}\ e_2+C_{1 3}^{3}\ e_3 \\{}
 [e_2,e_3 ]_{14} =C_{2 3}^{1}\ e_1\\
\end{array}$
\smallskip

$\begin{array}[t]{l}
 [e_1,e_2 ]_{15} =\frac{a\ C_{2 3}^{3}}{b}\ e_1+C_{1 2}^{2}\ e_2+\frac{C_{1 2}^{2}C_{2 3}^{3}}{C_{2 3}^{2}}\ e_3\\{}
 [e_1,e_3 ]_{15} = -\frac{((a-b)\ C_{1 2}^{2}-b \ C_{1 3}^{3})C_{2 3}^{2}}{b\ C_{1 2}^{2}}\ e_1 + C_{1 3}^{2}\ e_2+C_{1 3}^{3}\ e_3 \\{}
 [e_2,e_3 ]_{15} =\frac{a\ C_{2 3}^{2} C_{2 3}^{3}}{b\ C_{1 2}^{2}}\ e_1+C_{2 3}^{2}\ e _2+C_{2 3}^{3}\ e_3\\
\end{array}$
\smallskip

$\begin{array}[t]{l}
 [e_1,e_2 ]_{16} =\frac{a\ C_{2 3}^{3}}{b}\ e_1+\frac{C_{1 3}^{3}C_{2 3}^{3}}{C_{1 3}^{1}}\ e_3\\{}
 [e_1,e_3 ]_{16} = C_{1 3}^{1} e_1 + C_{1 3}^{2}\ e_2 +C_{1 3}^{3}\ e_3 \\{}
 [e_2,e_3 ]_{16} =\frac{a\ C_{1 3}^{2}C_{2 3}^{3}}{b\ C_{1 3}^{3}}\ e_1+C_{2 3}^{3}\ e_3\\
\end{array}$
\smallskip

$\begin{array}[t]{l}
 [e_1,e_2 ]_{17} =-\frac{ C_{1 3}^{1} C_{2 3}^{3}}{C_{2 3}^{2}}\ e_1-C_{1 3}^{3}\ e_2-\frac{C_{1 3}^{3}C_{2 3}^{3}}{C_{2 3}^{2}}\ e_3\\{}
 [e_1,e_3 ]_{17} = C_{1 3}^{3}\ e_1 +\frac{C_{1 3}^{3}C_{2 3}^{2}}{C_{2 3}^{3}}\ e_2+C_{1 3}^{3}\ e_3 \\{}
 [e_2,e_3 ]_{17} =C_{2 3}^{1}\ e_1+C_{2 3}^{2}\ e _2+C_{2 3}^{3}\ e_3\\
\end{array}$
\smallskip

$\begin{array}[t]{l}
 [e_1,e_2 ]_{18 }= C_{1 2}^{1}\ e_1 + C_{1 2}^{2}\ e_2 + C_{1 2}^{3}\ e_3 \\{}
 [e_1,e_3 ]_{18} =-\frac{a\ C_{2 3}^{2}}{b}\ e_1-\frac{C_{1 2}^{2}C_{2 3}^{2}}{C_{1 2}^{1}}\ e_2\\{}
 [e_2,e_3 ]_{18}=\frac{a\ C_{1 2}^{1}C_{2 3}^{2}}{b\ C_{1 2}^{2}}\ e_1+C_{2 3}^{2}\ e_2\\
\end{array}$
\smallskip

$\begin{array}[t]{l}
 [e_1,e_2 ]_{19} = -\frac{(b\ C_{1 2}^{2}+(b-a)\ C_{1 3}^{3})C_{2 3}^{3}}{b\ C_{1 3}^{3}}\ e_1 + C_{1 2}^{2}\ e_2 + C_{1 2}^{3}\ e_3 \\{}
 [e_1,e_3 ]_{19} =-\frac{a\ C_{2 3}^{2}}{b}\ e_1+\frac{C_{1 3}^{3}C_{2 3}^{2}}{C_{2 3}^{3}}\ e_2+C_{1 3}^{3}\ e_3 \\{}
 [e_2,e_3 ]_{19}=\frac{a\ C_{2 3}^{2}C_{2 3}^{3}}{b\ C_{1 3}^{3}}\ e_1+C_{2 3}^{2}\ e_2+C_{2 3}^{3}\ e_3\\
\end{array}$
\smallskip

$\begin{array}[t]{l}
 [e_1,e_2 ]_{20} = C_{1 2}^{2} \ e_1 + C_{1 2}^{2}\ e_2 +
 \frac{-a\ C_{1 3}^{3} C_{2 3}^{3}+b\ (C_{1 2}^{1}C_{1 3}^{1}+C_{1 2}^{2}+C_{1 3}^{3})C_{2 3}^{3}}{b\ C_{1 3}^{1}+a\ C_{2 3}^{2}}\ e_3 \\{}
 [e_1,e_3 ]_{20} = C_{1 3}^{1}\ e_1+\frac{b\ C_{1 2}^{2}C_{1 3}^{1}+ a\ C_{1 2}^{2}C_{2 3}^{2}+b\ C_{1 2}^{2}C_{2 3}^{2}-b\ C_{1 3}^{3}C_{2 3}^{2}}
 {b\ C_{1 2}^{1}+a\ C_{2 3}^{3}}\ e_2+C_{1 3}^{3}\ e_3 \\{}
 [e_2,e_3 ]_{20}=\frac{a\ C_{1 2}^{1}C_{2 3}^{2}+C_{1 3}^{1}C_{2 3}^{3}}{b\ (C_{1 2}^{2}+C_{1 3}^{3})}\ e_1+C_{2 3}^{2}\ e_2+C_{2 3}^{3}\ e_3\\
\end{array}$

\begin{rem}
The Lie algebras correspond in this list to the choice $a=b=1$.
\end{rem}


\end{document}